\documentclass[12pt, a4paper, reqno]{amsart}

\usepackage{graphicx}
\usepackage{enumerate}
\usepackage{amssymb}

\newtheorem{lemma}{Lemma}[section]
\newtheorem{corollary}[lemma]{Corollary}
\newtheorem{theorem}[lemma]{Theorem}
\newtheorem{proposition}[lemma]{Proposition}
\newtheorem{remark}[lemma]{Remark}

\newcommand{\Q}{\mathbb{Q}}

\author[A. Gasull]{A. Gasull$^*$}
\address{Dept. de Matem\`{a}tiques. Universitat Aut\`{o}noma de Barcelona. Edifici C. 08193 Bellaterra, Barcelona. Spain} 
\email{gasull@mat.uab.cat}

\author{H. Giacomini}
\address{Laboratoire de Math\'{e}matiques et Physique Th\'{e}orique. UMR 6083-CNRS- Facult\'{e} des
Sciences et Techniques. Universit\'{e} de Tours. 37200 Tours. France}
\email{Hector.Giacomini@lmpt.univ-tours.fr}

\author[J. Torregrosa]{J. Torregrosa$^*$}
\address{Dept. de Matem\`{a}tiques  Universitat Aut\`{o}noma de Barcelona Edifici C. 08193 Bellaterra, Barcelona. Spain} 
\email{torre@mat.uab.cat}

\thanks{$^*$Partially supported by a MCYT/FEDER grant number MTM2008-03437 and by a CIRIT
grant number 2009SGR 410.}


\begin{document}

\title[Homoclinic and heteroclinic connections]
{Some results on homoclinic and heteroclinic connections in planar
systems}

\begin{abstract}
Consider  a family of planar systems depending on two parameters
$(n,b)$ and having at most one limit cycle. Assume that the limit
cycle disappears at some homoclinic (or heteroclinic) connection
when $\Phi(n,b)=0.$  We present a method that allows to obtain a
sequence of explicit algebraic lower and upper bounds for the
bifurcation set ${\Phi(n,b)=0}.$ The method is applied to two
quadratic families, one of them is  the well-known Bogdanov-Takens
system. One of the results that we obtain for this system is the
bifurcation curve for small values of $n$, given by $b=\frac5 7
 n^{1/2}+\frac{72}{2401}n- \frac{30024}{45294865}n^{3/2}-
\frac{2352961656}{11108339166925} n^2+O(n^{5/2})$. We obtain the new
three terms from purely algebraic calculations, without eva\-luating
Melnikov functions.
\end{abstract}

\maketitle

\section{Introduction}

Consider a smooth  family of planar differential equations $(\dot
x,\dot y)=(P(x,y;n,b),\\Q(x,y;n,b)),$ $(n,b)\in \mathbb{R}^2,$  for which the
existence of at most one limit cycle is already known and moreover
all the bifurcations occurring in the family  are well understood.
For this family we could say that the {\it Qualitative theory of
ordinary differential equations} has achieved all its goals and all the job is done.
Nevertheless, from an analytic and quantitative point of view there
remains a crucial question: to determine the bifurcation curves in the
bifurcation diagram of the family.

Some  of these bifurcation curves are not difficult
to find. This is the case for instance of the curves that control the
changes in the behavior of the flow near the critical points.
These bifurcation curves correspond to local phenomena and so, in
principle, are easier to study. A paradigmatic example of this kind
is the bifurcation curve associated to the birth of a limit
cycle due to an Andronov-Hopf bifurcation. The
curves governing global phenomena are in general much more difficult
to determine. This is, for instance,  the case of the appearance of
homoclinic or heteroclinic connections.

We propose a method that, in some cases, allows to obtain explicit
lower and upper bounds in the parameter plane for the location of
the bifurcation curves associated to these global phenomena.   As
illustration,  we will apply it to  two 2-parametric families of
quadratic systems having at most one limit cycle. Since our examples
are polynomial, we will consider their phase portraits on the
Poincar\'{e} sphere,  see for instance \cite[Chap.~5]{DLA} or
\cite[Chap.~3.10]{Per2001},  and the homoclinic or heteroclinic
connections will also take into account the critical points at
infinity.

We start by giving a description of the method  for finding the relation between the
parameters $(n,b)$ for which this homoclinic or heteroclinic
connection exists. As a first step we consider an algebraic curve
$\mathcal {C}=\{(x,y)\in\mathbb{R}^2\,:\,C(x,y)=0\},$ as general as
possible, but having some of the geometrical properties of the phase
portrait of the system. For instance, it has to pass through
the two critical points of the vector field that are connected by
the heteroclinic orbit; if the critical point is a saddle it has to
be tangent to one of its separatrices; if the critical point is at
infinity, it has to have a branch going to infinity in the direction
corresponding to this point, etc.

The second step consists in using all the free parameters in
$\mathcal C$ to impose that, near the critical points, the curve be
as close as possible to the separatrix connecting orbit. In most
cases this can be done by finding a local analytical expression of
this separatrix in the neighborhood of each critical point. There
are several possibilities to impose that $\mathcal C$ be as close as
possible to the separatrix. If the separatrix is a heteroclinic
orbit we can fix all the free parameters of $\mathcal C$ imposing
this condition only in the neighborhood of one of the critical
points or in the neigbourhood of both critical points at the same
time. If the separatrix is a homoclinic orbit we can impose the
condition only on one of the branches of the orbit that arrives to
the critical point or on the two at the same time. In the two
families of quadratic systems studied in the next sections all these
possibilities will be explored.

Once all the parameters of the curve $\mathcal C$ are fixed, we
impose one more level of closeness between both objects (the
connecting orbit and the algebraic curve) by forcing a relation
among the parameters $(n,b)\in \mathbb{R}^2$. This procedure gives a
curve $\mathcal B =\{(n,b)\in \mathbb{R}^2 : \Phi(n,b)=0\}$ in the bifurcation space. In fact, up to this point the method proposed is essentially the one followed
in~\cite{BFG} to obtain analytic approximations to separatrices for some
two-point boundary value problems.

The last step is our main contribution and is the one that leads to
explicit algebraic upper and lower bounds of the bifurcation curve
associated to the existence of connections between the two critical
points. This  is the most difficult and computationally involved
part of the method. We have to  prove that there are parameters
$(n,b)$ near the curve $\mathcal B$ such that the corresponding
curves $\mathcal C$ are without contact for the flow of the vector
field and moreover that there are values for which the flow crosses
the curve in one sense and values for which the crossing is in the
opposite direction. This will be clearer in our applications.

Notice that the aim  of our work is in the spirit    of what Coppel
proposed in his well-known paper \cite{C}: ``Ideally one might hope
to characterize the phase portraits of quadratic systems by means of
algebraic inequalities on the coefficients", taking into account the
results of \cite{DF} where the authors proved that there are
bifurcation curves in quadratic systems which are neither algebraic
nor analytic. Since there is no hope to find analytic or algebraic
expressions of the bifurcation curves, we try to sandwich them
between two algebraic curves.

We present below our main results. The first one deals with a family
of quadratic systems already studied in \cite{Ro} and the second one
with the Bogdanov-Takens system. For both systems there is at most
one (hyperbolic) limit cycle. The value $b=b^*(n)$ corresponds, for
each $n$, to the value where this limit cycle, which is born in an
Andronov-Hopf bifurcation, disappears giving rise to a heteroclinic
or homoclinic connection, respectively, see more details in the next
sections.

\begin{theorem}\label{main1}
Let $b=b^*(n)$ be the bifurcation curve corresponding to the
heteroclinic connection for the quadratic system
\begin{equation}\label{sisbis}
\left\{
\begin{array}{ccl}
\dot x&=& P(x,y)=y, \\
\dot y&=&Q(x,y)=-x+by+xy-ny^2,
\end{array}
\right.
\end{equation}
with $b>0,$ $n>0$. Then:
\begin{enumerate}[(i)]

\item For $n>14,$
\[
\left|\frac1{4n} -b^*(n)\right|<\frac 1{8n^3}.
\]

\item  For $n>5,$
\[
\left|\frac1{4n}-\frac1{64n^3} -b^*(n)\right|<\frac 1{2n^5}.
\]
\item For $n>26,$
\[
\left|\frac1{4n}-\frac1{64n^3}-\frac{5}{512n^5} -b^*(n)\right|<\frac
1{8n^7}.
\]
\end{enumerate}
\end{theorem}

\begin{theorem}\label{mainbt} Let $b=b^*(n)$ be the bifurcation curve corresponding to
the saddle loop homoclinic connection for the Bogdanov-Takens system
\begin{equation}\label{sisbt}
\left\{
\begin{array}{ccl}
\dot x&=& y, \\
\dot y&=&-n+by+x^2+xy,
\end{array}
\right.
\end{equation}
with $n>0.$ Then $b^*(n)=H(\sqrt{n})$, where $H$ is an analytic
function. Moreover,
\begin{enumerate}[(i)]

\item For $0\le n \le 1/40,$
\begin{equation*}
\left|\frac 57  \sqrt n+ \frac{72}{2401}n-b^*(n)\right|<\frac{11}4
n^{5/4}.
\end{equation*}

\item For $n$ small enough,
\begin{align*}
b^*(n)&=\frac5 7 n^{1/2}+\frac{72}{2401}n-
\frac{30024}{45294865}n^{3/2}- \frac{2352961656}{11108339166925}
n^2+O(n^{5/2}).
\end{align*}
\end{enumerate}

\end{theorem}

The fact that $b^*(n)=H(\sqrt{n}),$  where $H$ is an analytic
function, is an easy consequence of some results of~\cite{Perko92},
see Corollary~\ref{ccc}. Our main contributions are items (i) and
(ii). Notice that  they improve the well-known local knowledge of
the function $b=b^*(n)=\frac5 7\sqrt{n}+O(n),$ proved for instance
in \cite{GH} by using a blowing-up process and scaling of the time
and the computation of a Melnikov function. In Section~\ref{s3.3} we
compare them with other results of Perko obtained in the mentioned
paper~\cite{Perko92}. Finally, in Subsection~\ref{seckorea} we prove
that
\[
b^*(1/4)\in\left({\textstyle\frac{951225059}{2609347034},
 \frac{258052528}{707875165}}\right)\simeq
 (0.364545247,0.3645452486).\]
We use this information to improve the results of~\cite{korea}. This
subsection also serves  to show how to use our approach when the
problem that we face depends only on one parameter.

At the best of our knowledge, it is the first time that this type of
results are obtained for global bifurcation problems.

\section{The first family of quadratic systems}

\begin{figure}[h]
\begin{center}
\includegraphics[scale=0.35]{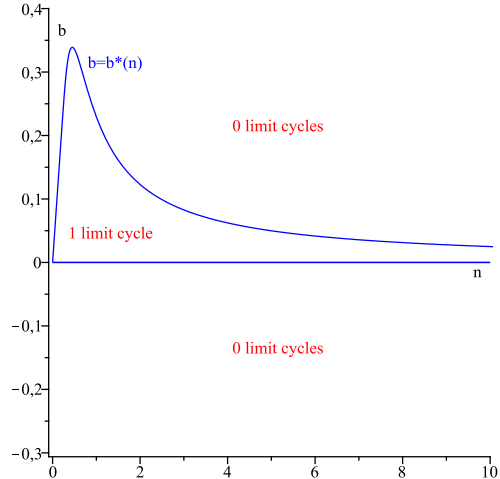}
\end{center}
\caption{Limit cycles of system~\eqref{sisbis} for $n>0.$ The curve
$b=b^*(n)$ is obtained numerically.}\label{fig1}
\end{figure}

Our starting point is the following theorem:

\begin{theorem}\label{known} Consider the system~\eqref{sisbis},
\begin{equation*}
\left\{
\begin{array}{ccl}
\dot x&=& P(x,y)=y, \\
\dot y&=&Q(x,y)=-x+by+xy-ny^2,
\end{array}
\right.
\end{equation*}
with $(n,b)\in\mathbb{R}^2.$ Then there exits a function $b=b^*(n)$
such that $b^*(-n)=-b^*(n),$ $n b^*(n)>0$ for $n\ne0$ for which:
\begin{enumerate}[(i)]

\item  The system has exactly one hyperbolic limit cycle if either
$n>0$ and $0<b<b^*(n)$ or $n<0$ and $b^*(n)<b<0.$

\item  Otherwise the system has no limit cycles.
\end{enumerate}
See Figure~\ref{fig1} for a numeric plot of the function $b=b^*(n)$.
\end{theorem}

Essentially its proof follows from \cite{Ro} and the theory of
rotated vector fields \cite{Duff,Perko}. See also \cite{CGL} and
\cite{Ro}. The property $b^*(-n)=-b^*(n)$ is a consequence of the
fact that the transformation $(x,y,t)\rightarrow (-x,y,-t)$ changes
the sign of $b$ and $n$ in system~\eqref{sisbis}.

\begin{figure}
\begin{center}

\begin{tabular}{cccc}
\includegraphics[height=3cm]{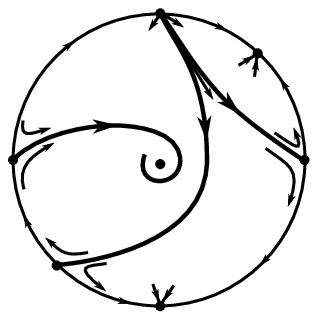}&\includegraphics[height=3cm]{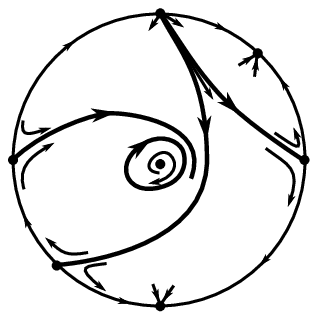}&
\includegraphics[height=3cm]{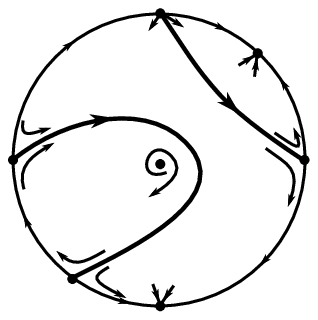}&\includegraphics[height=3cm]{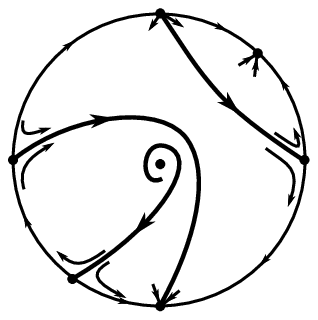}\\
$b\leq0$& $0<b<b^*(n)$&$b=b^*(n)$& $b^*(n)<b<n$\\
\end{tabular}
\end{center}
\caption{Phase portraits of system~\eqref{sisbis} for $n>0$ and
$b<n$. }\label{fig2}
\end{figure}

The phase portraits of system~\eqref{sisbis} in the Poincar\'{e} sphere
appear in \cite{CGL}. Indeed the phase portrait corresponding to
$b=b^*(n),$ for $n\ne0,$ corresponds to the one having a
heteroclinic connection of the two infinite semi-hyperbolic critical
points which are on the directions $y=0$ and $x-ny=0,$ see
Figure~\ref{fig2}.  From the above considerations we already know
that we can restrict our attention to the region $n>0$ and $b>0.$
The origin is the only finite critical point of the system. Moreover
the line $y-1=0$ is transversal to the flow associated to this
system. A key result for knowing, once $n$ is fixed, if some value
$b$ is smaller or bigger than $b^*(n)$ is the following lemma.

\begin{figure}[h]
\begin{center}
\begin{tabular}{ccc}
\includegraphics[height=3cm]{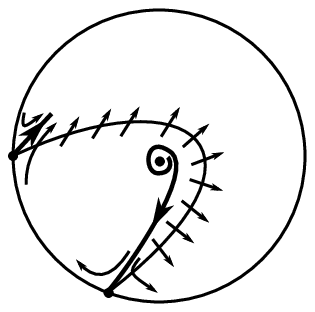}&\phantom{xxxxxxx} &\includegraphics[height=3cm]{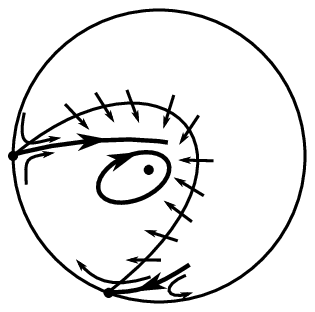}\\
$M<0$&\phantom{xxxxxxxx} & $M>0$\\
\end{tabular}

\end{center}
\caption{Flow on the curve $\phi(x,y)=0$.}\label{fig3}
\end{figure}

\begin{lemma}\label{contacto} Let $\phi(x,y)=q(y)x-p(y)=0,$ with $p(y)$ and $q(y)$ polynomials,
be an algebraic curve having a non-singular branch connecting the
two critical points at infinity of system~\eqref{sisbis}, which
correspond with the directions $y=0$ and $x-ny=0.$ Assume also that
$\phi(0,0)>0$. If for some $n>0$ and $b>0$ it holds that
\[
M(x,y):=\left. \frac{\partial \phi(x,y)}{\partial
x}P(x,y)+\frac{\partial \phi(x,y)} {\partial
y}Q(x,y)\right|_{\{\phi(x,y)=0\}}\ge0 \quad (\mbox {resp.} \le0)
\]
and the zeroes of $M$ are isolated, then $b<b^*(n)$ (resp.
$b>b^*(n)$), see Figure~\ref{fig3}.
\end{lemma}

\begin{proof} Consider for instance the case  where $M$ is
positive or zero. The inequality $M(x,y)\ge0$ implies that the
connected component of the curve $\{\phi(x,y)=0\}$ described in the
statement is crossed by the flow of the  vector field $(P,Q)$ from
the region $\{\phi(x,y)<0\}$ to the  region $\{\phi(x,y)>0\}.$ Since
for $b>0$ the origin is a repellor, by applying the
Poincar\'{e}-Bendixson Theorem on the sphere to the positive invariant
region formed by the curve and a piece of the equator (which is
invariant by the flow), we can prove the existence of at least one
periodic orbit of system~\eqref{sisbis}. By using
Theorem~\ref{known} we get that $b<b^*(n),$ as we wanted to prove.
The other case follows similarly.
\end{proof}

The main goal of this section is to get good lower and upper bounds
for the function $b=b^*(n)$.

\subsection{A first approach to the bifurcation
curve}

\begin{figure}[h]
\begin{center}
\includegraphics[scale=0.35]{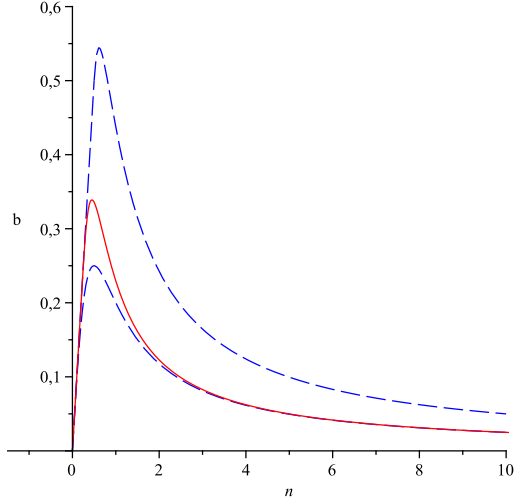}
\end{center}
\caption{Observe how the bifurcation curve $b=b^*(n)$ of
system~\eqref{sisbis}, obtained numerically, lies
  between the two dashed curves given in Proposition~\ref{prop1}.}\label{fig4}
\end{figure}

Our first result takes advantage of a direct application of
Lemma~\ref{contacto} by taking the curves $\phi(x,y)=0$ as
hyperbolas. As can be seen in Figure~\ref{fig4}, the curve $b=b^*(n)$ is quite well
delimited by this simple
approach.

\begin{proposition}\label{prop1}
Let $b=b^*(n)$ be the bifurcation curve for the quadratic
system~\eqref{sisbis}, with $n>0$ and $b>0$. Then, for all $n>0$:

\[\frac{n}{4n^2+1}<b^*(n)<
\begin{cases}
n,\quad  &\mbox{when}\quad n\le 1/2,\\
\dfrac{8n^2-1}{16n^3},\quad  &\mbox{when}\quad n\ge 1/2.
\end{cases}
\]
\end{proposition}

\begin{proof}
Consider $\phi(x,y)=a_0+a_1x+a_2y+xy-ny^2.$ Clearly, inside this
family of hyperbolas there are many of them having the geometrical
restrictions described in Lemma~\ref{contacto}. By computing
$M(x,y),$ since on $\phi(x,y)=0,$ we can isolate $x$ in terms of
$y$, we get that $M$  depends only on $y$. Moreover its numerator is
a polynomial of degree 4 in $y$ and its denominator is a polynomial of degree 2 which is a perfect square. There are several ways of forcing
it to have a given sign in the region where the connected component
of $\phi(x,y)=0$ of our interest lies. For instance, by taking
$n>1/2,$
\[a_0=\frac{1+4n^2}{8 n^3},\,a_1=-1,\,
a_2=-\frac{1+8n^2}{16n^3}\quad \mbox{and}\quad
b=\frac{8n^2-1}{16n^3},\] we force this polynomial to be $(y-1)$
times a polynomial of degree two which is a perfect square.
Concretely we get that
\[M(x(y),y)=\frac{(4n^2+1)(4n^2y-1)^2}{128 n^6(y-1)}<0\] in the region
$\{y<1\}$ and that the hyperbola connects the two infinite critical
points, as shown in Figure~\ref{fig3}. By using Lemma~\ref{contacto}
we get that for $n>1/2,$ ${(8n^2-1)}/{(16n^3)}>b^*(n),$ as we wanted
to prove.

By imposing that the polynomial in the numerator of $M(x(y),y)$ be
the product of $y^2$ and another polynomial of degree two, which is
also a perfect square, we get
\[
\begin{aligned}
&a_0=\frac{(2n^2+1)^2(8n^2+1)^2}{2n(4n^2+1)^4},\, &a_1=-\frac{(2n^2+1)(8n^2+1)}{(4n^2+1)^2},\\
&a_2=-\frac{(2n^2+1)(8n^2+1)}{2n(4n^2+1)^2},\,  &b=\frac n{4n^2+1}.
\end{aligned}
\]
Moreover $M(x(y),y)=y^2R_1^2(y)/S_1^2(y)>0,$ where $R_1$ and $S_1$
are polynomials of degree one in $y,$  and $S_1(y)\ne0$ in the
region $\{y<1\}.$ Again by Lemma~\ref{contacto} we get that for all
$n>0,$ $n/(4n^2+1)<b^*(n).$

To end the proof it only remains to show that for $0<n<1/2,$
$b^*(n)<n.$ Indeed this is true for all $n>0.$ Notice that the
straight line $\{y=1\}$ can not be cut by the periodic orbits of the
system. Since,  when $b-n\ge0,$
\[\mbox{div}\left(\frac{y}{1-y},\frac{-x+by+xy-ny^2}{1-y} \right)=
\frac{b-2ny+ny^2}{(1-y)^2}\ge0,\]  by applying the Bendixson-Dulac
criterion we get that the system has no limit cycles. Thus
$b^*(n)<n,$ as we wanted to prove.
\end{proof}

\begin{remark}
By imposing that the polynomial of degree four appearing in the
above proof be a multiple of $y^4$ we get another quite satisfactory
lower bound for $b=b^*(n)$ given by the function
$b=(4n^2-1)/(16n^3)$. Notice that $(4n^2-1)/(16n^3)<n/(4n^2+1).$
\end{remark}

\subsection{Proof of Theorem~\ref{main1}}

For simplicity, we only give all the details of the proof of the
weaker result, which is the one given in item (i), which asserts
that for $n>14,$
\[
\left|\frac1{4n} -b^*(n)\right|<\frac 1{8n^3}.
\]
Nevertheless its proof includes all the ingredients for  the proof
of the other results. To see the main difference for proving the
other two statements see the footnote of next page. Before starting
the proof we want to comment that although it seems natural to
improve item (iii) of the theorem going further in our computations,
our results are on the limit of what we can do with the algebraic
manipulator that we use (Maple) and the capacity of our
computers\footnote{We have used a computer with  32Gb of RAM and
2.4GHz.}.

To prove item (i) we will apply in this case the general method
described in the previous section. Recall that when $b=b^*(n)$ there
is a heteroclinic connection between two critical points at infinity
of the Poincar\'{e} compactification of system~\eqref{sisbis}, see
Figure~\ref{fig2}. Recall also that the directions corresponding to
these critical points are given by the lines
$R(x,y):=xQ_2(x,y)-yP_2(x,y)=0,$ where $P_2$ and $Q_2$ are the
quadratic homogeneous parts of $P$ and $Q,$ respectively. In our
case $R(x,y)=xy(x-ny).$ The behaviour near infinity of the
separatrices of these critical points can be easily obtained. For
instance the one corresponding to the direction $x-ny=0,$ from the
center manifold theorem, can be written as
\[
x=\Psi(y)=ny+\psi\left(\frac1y\right),\] with  $\psi$ being an analytic
function at zero of the form
\[\psi(w)=\frac{
1+{n}^{2}-nb}{n} \left( 1+w+\frac{n^2+1}{n^2}w^2+\frac {( n^2+2)(
n^2+1)}{n^4} w^3+O(w^4)\right).
\]
Indeed we will use this expansion at  much higher order.

Now we search an algebraic curve of the form\footnote{\label{foot2}To prove
statements  (ii) and (iii), the only difference is that we take a
rational function with the same structure but having higher degrees
in both the numerator and the denominator. Concretely, having
respectively degrees 6 and 5 for (ii), and 8 and 7 for (iii).}
\begin{equation*}\label{f4}
x=f(y):=\frac{p_u(y)}{q_v(y)}={\frac
{u_{{0}}+u_{{1}}y+u_{{2}}{y}^{2}+u_{{3}}{y}^{3}+u_{{4}}{y}^{4}+n{y}^{5}}{ \left( y-1 \right)  \left(
v_{ {0}}+v_{{1}}y+v_{{2}}{y}^{2}+{y}^{3} \right)
}}=ny+\sum_{k=0}^\infty f_k y^{-k},
\end{equation*}
which is as close as possible to the searched separatrix. This is
done by imposing that $\Psi(y)$ and $f(y)$ coincide at the highest
order possible at infinity. Concretely, the eight free constants in
$f(y)$ are used to vanish the coefficients of $y^{-k},
k=0,1,\ldots,7$ of $\Psi(y)-f(y),$ obtaining the asymptotic
expansion
\[
\Psi(y)-f(y)=\frac{C_8(n,b)}{y^8}+\frac{C_9(n,b)}{y^9}+\cdots,
\]
where
\[
C_8(n,b)=\frac { \left( -1-{n}^{2}+nb \right)N_8(n,b)}{D_8(n,b)},
\]
being
\[
\begin{aligned}
N_8(n,b)=&4{n}^{19}b-23{n}^{18}{b}^{2}+70{n}^{17}{b}^{3}-87{n}^{16}{b}^{4}+36{n}^{15}{b}^{5}-{n}^{18}-231{n}^{17}b\\
&+475{n}^{16}{b}^{2}-143{n}^{15}{b}^{3}-100{n}^{14}{b}^{4}+59{n}^{16}-2488{n}^{15}b+3668{n}^{14}{b}^{2}\\
&-886{n}^{13}{b}^{3}-232{n}^{12}{b}^{4}+567{n}^{14}-11151{n}^{13}b+10886{n}^{12}{b}^{2}-1284{n}^{11}{b}^{3}\\
& +2561{n}^{12}-28350{n}^{11}b+16332{n}^{10}{b}^{2}-272{n}^{9}{b}^{3}+7718{n}^{10}-44884{n}^{9}b\\
&+12264{n}^{8}{b}^{2}+17160{n}^{8}-44920{n}^{7}b+3648{n}^{6}{b}^{2}+27000{n}^{6}-26112{n}^{5}b\\
&+27024{n}^{4}-6624{n}^{3}b+14976{n}^{2}+3456
\end{aligned}
\]
and
\[
\begin{aligned}
D_8(n,b)=&-9{n}^{9}b+9{n}^{8}{b}^{2}-33{n}^{7}b+8{n}^{
6}{b}^{2}+2{n}^{8}+16{n}^{6}\\&-52{n}^{5}b+50{n}^{4}-28{n}^{3}b+60{n}^{2}+24
.
\end{aligned}
\]
Therefore, given $n=\bar n,$  some of the values of $b$ satisfying
$N_8( \bar n,b)=0$ could be a good candidate to approximate
$b^*(\bar n).$ It is not difficult to prove that one of the
connected components of $N_8(n,b)=0$  has an asymptotic expansion at
infinity of the form $b=\frac1{4n}+\frac{3}{64n^3}+\cdots.$ From
this expansion, and taking into account the results of
Proposition~\ref{prop1}, we decide to fix
\[
b=\frac{1}{4n}+\frac{\alpha}{n^3},
\]
for some $\alpha$ to be determined afterwards and the
$u_i=u_i(n,\alpha)$ and $v_i=v_i(n,\alpha)$ obtained above.
Taking $\phi_{n,\alpha}(x,y)=q_v(y) x-p_u(y)$ in
Lemma~\ref{contacto} we get that
\begin{equation}\label{final}
M_{n,\alpha}(f(y),y)=t(n,\alpha)
\frac{\sum_{i=0}^{4}r_i(n,\alpha)z^i}{z\left(\sum_{i=0}^{3}s_i(n,\alpha)z^i\right)^2},
\quad\mbox{where}\quad z=1-y,
\end{equation}
and  $f(y)=p_u(y)/q_v(y)$. In the above expression all the
functions $r_i(n,\alpha)$ are polynomials  in $n$ of degree 40, and all
$s_i(n,\alpha)$ are  also polynomials  in $n.$ The function $t(n,\alpha)$ is a
rational function, which is always negative, and the coefficients of
$n^{40}$ for each $r_i(n,\alpha),$ for $i=0,\dots,4$ are, respectively,
\[
-\frac{5}{256}-\frac{\alpha}{4},\,-\frac{3}{64}-\alpha,\,-\frac{3}{128}
-\frac{3\alpha}{2},\,\frac{1}{64}-\alpha,\,\frac{3}{256}-\frac{\alpha}{4}.
\]

When $\alpha=1/8$ (resp. $-1/8$) all the expressions of the above list
are negative (resp. positive). Hence  we can assure that for this
value of $\alpha$ there exist $n_0^{\pm}$ such that for all $z>0$  and all
$n>n_0^+$ (resp. $n>n_0^-$) the numerator of \eqref{final} is
positive (resp. negative). Similarly it can be proved that the
denominator never vanishes on the same region. Indeed it can be seen
that the biggest value of $n$ which vanishes some of the eighteen
functions, $r_i(n,\pm1/8),i=0,\dots,4,$ and
$s_i(n,\pm1/8),i=0,\dots,3$ corresponds to a zero of a factor of
$r_4(n,1/8).$ This factor is
\[
\begin{aligned}
&2560{n}^{16}-446720{n}^{14}-2294624{n}^{12}+1116256 {n}^{10}+33662656 {n}^{8}\\
&+98872176{n}^{6}+131940378{n}^{4}+85765842{n}^{2}+21990713,
\end{aligned}
\]
and its biggest zero  is approximately $13.397.$ In short, by using
for instance the Sturm algorithm, it can be proved that for $\alpha=1/8$
(resp. $\alpha=-1/8$), $n>13.4$ and $z>0$ it holds that the expression in
\eqref{final} is  positive (resp. negative). By using
Lemma~\ref{contacto} the theorem follows.

\begin{remark} It is clear that our method to choose the
curve $\phi(x,y)=0,$ which is one of the key points of our approach,
that consists in imposing that it coincides as much as possible with
one of the separatrices, can be done by using the separatrix of the
other infinite critical point. Also it can be imposed that the curve
$\phi(x,y)=0$ coincides  in both extremes, i.e. simultaneously as
much as possible, with both separatrices. Although we obtain also
other bounds for $b=b^*(n),$ the ones given in Theorem~\ref{main1}
are, at least for $n$ big, better than the ones obtained by
using these slightly different approaches.
\end{remark}

\section{The Bogdanov-Takens system}\label{btseccion}

\begin{figure}[h]
\begin{center}
\includegraphics[scale=0.35]{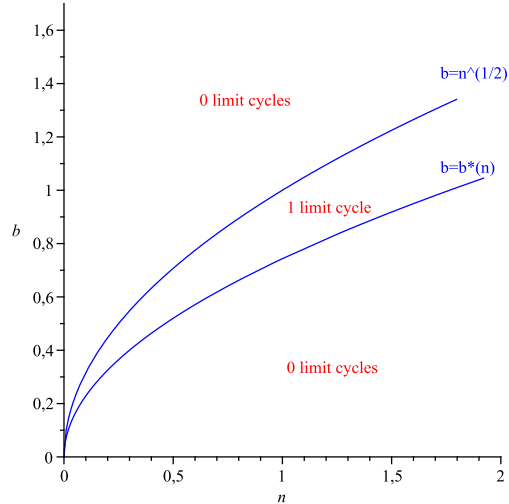}
\end{center}
\caption{Limit cycles of the Bogdanov-Takens system. The curve
$b=b^*(n)$ is obtained numerically.}\label{fig5}
\end{figure}

For the Bogdanov-Takens  system all its bifurcation diagram,
including the number of limit cycles is known, see for instance
\cite{CLW, GH, LRW, Perko92}. We summarize these results about
limit cycles in the next theorem:

\begin{theorem}\label{teobt} Consider the system~\eqref{sisbt},
\begin{equation*}
\left\{
\begin{array}{ccl}
\dot x&=& y, \\
\dot y&=&-n+by+x^2+xy,
\end{array}
\right.
\end{equation*}
with $n$ and $b$ real numbers. Then it has  a limit cycle (which is
unique and hyperbolic) if and only if
\[
n>0 \quad \mbox{and} \quad b^*(n)<b<\sqrt{n}.
\]
Here the function $b^*\hspace{-0.1cm}:\mathbb{R}^+\to \mathbb{R}$
gives the curve $b=b^*(n),$ $n\ge0$, where the phase portrait of
\eqref{sisbt} has a homoclinic loop, which is hyperbolic and
unstable. Moreover
\[
b^*(n)=\frac 5 7 \sqrt{n}+O(n),\quad n\sim 0.
\]
See Figure \ref{fig5} for a numeric plot of $b=b^*(n).$
\end{theorem}

\begin{figure}[h]
\begin{center}
\begin{tabular}{cccc}
\includegraphics[height=1.7cm]{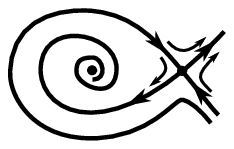}&\includegraphics[height=1.7cm]{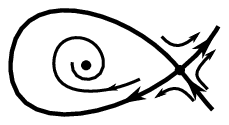}&
\includegraphics[height=1.7cm]{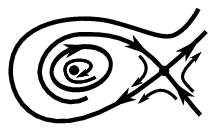}&\includegraphics[height=1.7cm]{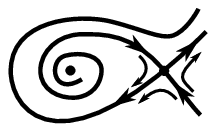}\\
$b< b^*(n)$& $b=b^*(n)$&$b^*(n)<b<\sqrt{n}$& $b\geq\sqrt{n}$\\
\end{tabular}
\end{center}
\caption{Saddle loop and Hopf   bifurcations for the Bogdanov-Takens
system.}\label{fig6}
\end{figure}

Similarly that  in the previous section our main goal is to  give
analytic information on the function $b=b^*(n)$ above. Following
\cite{CLW,GH}  we know that, fixed $n>0$, a hyperbolic unstable
limit cycle is created when $b=\sqrt n,$ via a Hopf bifurcation,
that this limit cycle increases size when $b$ decreases and
disappears for some value $b=b^*(n),$ for which the system presents
a homoclinic unstable saddle loop. This behavior of the limit cycle
is due to the fact that the system is a semi-complete family of
rotated vector fields with parameter $b$, see for instance
\cite{Duff, Perko}, see Figu\-re~\ref{fig6}. The uniqueness and
hyperbolicity of the limit cycle is proved for instance
in~\cite{LRW}.

\subsection{A first approach to $b^*(n)$}

When $n\ge0,$ instead of working with the expression~\eqref{sisbt},
by using the translation $x_1=x-\sqrt{n}$, $y_1=y$, after dropping
the subindexes, we obtain the equivalent system
\begin{equation}\label{bt}
\left\{
\begin{array}{ccl}
\dot x&=& y, \\
\dot y&=&2\sqrt{n}x+(b+\sqrt{n})y+x^2+xy.
\end{array}
\right.
\end{equation}

A simple first estimation of the function $b=b^*(n)$ is given in
next lemma:

\begin{lemma}\label{btprel} Let $b=b^*(n)$ be the value corresponding to a homoclinic saddle loop  for system~\eqref{sisbt}, then
\[
\max(-\sqrt n, \sqrt {n}-1)< b^*(n)< \sqrt{n}.
\]
\end{lemma}

\begin{proof} We work  with the  expression~\eqref{bt} of the Bogdanov-Takens system.
Recall that the limit cycle is born, via a Hopf bifurcation, when
$b=\sqrt n$ and $b$ decreases. Hence, by using the non-intersection
property of the limit cycles of a semi-complete family of rotated
vector fields it suffices to prove that, when either $b=\sqrt{n}-1$
or $b=-\sqrt n$, system \eqref{bt} has no limit cycles. When
$b=\sqrt{n}-1$, note that the line $x+y=0$ is invariant by the flow
of \eqref{bt}, because
\begin{equation}\label{ll}
\dfrac{\partial L(x,y)}{\partial x}y+\dfrac{\partial L(x,y)}
{\partial
y}\left(2\sqrt{n}x+(b+\sqrt{n})y+x^2+xy\right)=\left(2\sqrt{n}+x\right)
L(x,y),
\end{equation}
where $L(x,y)=x+y.$ In other words, $L(x,y)$ is an invariant
algebraic curve with cofactor $K(x,y):=2\sqrt{n}+x.$

Assume that, when $b=\sqrt{n}-1$, system~\eqref{bt} has a limit
cycle $\Gamma=\{(x(t),y(t))\},$ with period $T$. Since the
divergence of the system is $b+\sqrt{n}+x=2\sqrt{n}-1+x,$ the
characteristic exponent of $\Gamma$ is
\[
\kappa:=\int_0^T \left(2\sqrt{n}-1+x(t)\right) \,dt.\] On the other
hand notice that from~\eqref{ll} we get that
\[
\int_0^T \left(2\sqrt{n}+x(t)\right)
\,dt=\int_0^T\frac{\frac{d}{dt}L(x(t),y(t))}{L(x(t),y(t))}=0.
\]
Hence $\kappa=-T$ and so $\Gamma$ would be a hyperbolic stable limit
cycle. This fact is in contradiction with the known results, see for
instance Figure~\ref{fig6}. Hence the limit cycle does not exist in
this case.

To end the proof let us show that $b^*(n)>-\sqrt n.$ It is
well-known that the stability of a hyperbolic saddle loop is given
by the sign of the divergence at the saddle point, see \cite[Chap.
XI]{ALGM}. In our case, the divergence at the saddle point is
$b^*(n)+\sqrt{n}$. For the Bogdanov Takens system \eqref{sisbt},
when the loop exists it is hyperbolic and unstable. Hence it holds
that $b^*(n)+\sqrt n>0$, as we wanted to prove.
\end{proof}

\subsection{Proof  of Theorem~\ref{mainbt}.(i)}
We introduce the new parameter $m=n^{1/4}$ and we consider the
equivalent system given in~\eqref{bt},
\begin{equation}\label{btm}
\left\{
\begin{array}{ccl}
\dot x&=& y, \\
\dot y&=&2m^2x+(b+m^2)y+x^2+xy.
\end{array}
\right.
\end{equation}

We will prove our result in a constructive way. Indeed our first steps
could be skipped but we believe that they are useful to understand the
way we have obtained the approximations for $b^*(n).$

System \eqref{btm} has a focus or a node at the point $(-2m^2,0)$ and a
saddle point at the origin. The linear approximation to its
separatrices at the origin is given by the two lines:
\[
C_2(x,y):=-2 m^2 x^2 -(b+m^2)xy+y^2=0,
\]
or equivalently, $y= a_1x,$ where $a_1$  is one of the values
\begin{equation}\label{apm}
a_1^{\pm}:=\frac{b+m^2\pm\sqrt{(b+m^2)^2+8m^2}}{2}.
\end{equation}

As in the proof of Theorem~\ref{main1} we start by finding the first
terms in the Taylor expansion of the separatrices of a critical
point. In this occasion the critical point is the origin. In its
neighborhood we can write these separatrices as
\[y=\Psi(x):=\sum_{k=1}^\infty a_k x^k,\]
where $a_1$ is one of the values $a_1^{\pm}$, or in other words, $C_2
(1,a_1)=0.$ The next terms can be found recurrently by imposing that
\[
\left.\frac{\partial (y-\Psi(x))}{\partial x}y+\frac{\partial
(y-\Psi(x))}{\partial y}\left(
2m^2x+(b+m^2)y+x^2+xy\right)\right|_{y=\Psi(x)}\equiv0.\] We have,
for instance,
\begin{align*}
a_2&={\frac {-1-a_{{1}}}{b-3a_{{1}}+{m}^{2}}}, &a_3&={\frac {a_{{2}}
\left( 2a_{{2}}-1 \right) }{b-4a_{{1}}+{m}^{2}}},\\
a_4&={\frac {a_{{3}} \left( -1+5a_{{2}} \right)
}{b-5a_{{1}}+{m}^{2}}},& a_5&={\frac
{3{a_{{3}}}^{2}-a_{{4}}+6a_{{2}}a_{{4}}}{b-6a_{{1}}+{m}^ {2}}}.
\end{align*}
We remark  that $b-k a_1+m^2$ never vanishes for $k>1.$

Now we consider a cubic algebraic curve of the form
\begin{equation}\label{cubica}
C(x,y):=C_2 (x,y)+c_{3,0}x^3+c_{2,1}x^2y+c_{1,2}xy^2+c_{0,3}y^3=0
\end{equation}
and we search the four free coefficients by imposing that $\{C=0\}$
be as close as possible to the separatrices of the saddle point. In
particular notice that the quadratic terms of $C$ imply that this
curve is tangent to both separatrices. This can be done by imposing
that the function
\[
F(x):=C(x,\Psi(x))=\sum_{k=3}^\infty f_k(m,b) x^k,
\]
be as flat as possible at the origin. By choosing suitable $c_{i,j}$
 in terms of $b$ and $m$ we get that
\begin{equation*}
F(x)=f_7(m,b) x^7+O(x^8).
\end{equation*}
After some cumbersome computations, done with an algebraic
manipulator,  we obtain that all the solutions of the equation
$f_7(m,b)=0$ are contained in the algebraic curve

\begin{align*}
& (2+b)(b+1-m^2)(b^2+9m^2+2bm^2+m^4)^2\times\\&\big(3726
{m}^{20}+21951 {m}^{18}+25110 b{m}^{18}-671220 {m}^{16}+68040
 {b}^{2}{m}^{16}\\
&+109512 b{m}^{16}+169452 {b}^{2}{m}^{14}+87480 {b}
^{3}{m}^{14}-3005247 {m}^{14}\\
&-1882701 b{m}^{14}-1560330 {b}^{2}{m}^ {12} -40176
{b}^{3}{m}^{12}-5046048 b{m}^{12}\\
&+17286120 {m}^{12}+ 34020 {b}^{4}{m}^{12}-703539
{b}^{3}{m}^{10}-432702 {b}^{4}{m}^{10}\\
& -2360050 {m}^{10}-1021503 {b}^{2}{m}^{10}-39074153 b{m}^{10}-47628
 {b}^{5}{m}^{10}\\
&-848112 {b}^{3}{m}^{8}+31843300 {b}^{2}{m}^{8}+ 640000
{m}^{8}-530064 {b}^{5}{m}^{8}\\
&-1329192 {b}^{4}{m}^{8}-68040 {b}^{6}{m}^{8}+9121240
b{m}^{8}-1181763 {b}^{5}{m}^{6}\\
&-29160 {b}^{7 }{m}^{6}-1792000 b{m}^{6}-238788
{b}^{6}{m}^{6}-2892813 {b}^{4}{m}^ {6}\\
&-11432588 {b}^{2}{m}^{6}-18903846 {b}^{3}{m}^{6}-89544 {b}^{5}{m
}^{4}+1254400 {b}^{2}{m}^{4}\\
&+22032 {b}^{7}{m}^{4}+56862 {b}^{6}{m}^ {4}+4247140
{b}^{4}{m}^{4}+2430 {b}^{8}{m}^{4}\\
&+4429720 {b}^{3}{m}^{ 4}+974859 {b}^{6}{m}^{2}+5670
{b}^{9}{m}^{2}+243950 {b}^{4}{m}^{2}\\
&+54999 {b}^{8}{m}^{2}+307107 {b}^{7}{m}^{2}+863815 {b}^{5}{m}^{2}
+1296 {b}^{10}\\
&+13608 {b} ^{9}+42984 {b}^{8}+39000 {b}^{7}+11400 {b}^{6}\big)=0.
\end{align*}

The branches passing through the origin of the algebraic curve given above are
\begin{equation*}
b=\frac5 7
m^2+\frac{72}{2401}m^4\pm\frac{4428}{2401}\sqrt{2}\,m^5+O(m^6).
\end{equation*}
Hence to continue our study we decide to fix
\begin{equation}\label{bfix}
b=\frac5 7 m^2+\frac{72}{2401}m^4+ \alpha\, m^5,
\end{equation}
with $\alpha$ to be given afterwards.

With $b$ given  in \eqref{bfix} we will try to fix $\alpha$ in such a way
that the curve $C$ be without contact in the negative half-plane
$\{(x,y)\,:\, x<0\}$ where the saddle separatrices lie. To impose
this constrain we compute the resultant  of $C(x,y)$ and
\begin{align*}
D(x,y)&:=\frac{\partial C(x,y)}{\partial x}y+\frac{\partial
C(x,y)}{\partial y}\left( 2m^2x+(b+m^2)y+x^2+xy\right),
\end{align*}
with respect to $y.$ This a very huge task. After many
computations\footnote{Notice that this is a main difference between
the proof of this theorem and the one of Theorem~\ref{main1},
because there the algebraic curve considered is of the form
$x=f(y)$, and so the problem goes easily to a one variable problem,
while in this case we have to evaluate the resultant of two
polynomials.}, done again with an algebraic manipulator, we obtain
that
\[
P(x;\alpha,m):=
\textrm{Res}(C,D,y)=x^{10}m^{36}\left(r_0(\alpha,m)+r_1(\alpha,m)x+r_2(\alpha,m)x^2\right),
\]
where for $i=0,1,2,$ it holds that
\[r_i(\alpha,m)=r_i^0(\alpha,m)+r_i^1(\alpha,m)\Delta(\alpha,m),
\]
where
\begin{align*}
\Delta(\alpha,m):=&\big(5764801 {\alpha}^{2}{m}^{8}+345744
\alpha\,{m}^{7}+5184 {m}^{6}+19765032 \alpha\,{m}^ {5}\\&+592704
{m}^{4}+16941456 {m}^{2}+46118408\big) ^{1/2}
\end{align*}
and all the functions $r_i^j(\alpha,m)$  are polynomials in $\alpha$ and $m$.
Solving the quadratic equation
\[
p_2(x;\alpha,m):=r_0(\alpha,m)+r_1(\alpha,m)x+r_2(\alpha,m)x^2=0,
\]
we obtain the solutions
\begin{align*}
x=x^+(\alpha,m)&=\frac{-r_1(\alpha,m)+\sqrt{r_1^2(\alpha,m)-4 r_0(\alpha,m)r_2(\alpha,m)}}{2 r_2(\alpha,m)}\\
&=\left(-\frac{63}{25}+\frac{16807}{18450}\sqrt{2}\,\alpha+\frac{282475249}{326786400}\,\alpha^2
\right)m^2+O(m^3),\\
x=x^-(\alpha,m)&=\frac{-r_1(\alpha,m)-\sqrt{r_1^2(\alpha,m)-4 r_0(\alpha,m)r_2(\alpha,m)}}{2 r_2(\alpha,m)}\\
&=\frac{247107}{304580}+O(m).
\end{align*}

Notice that the two roots of the coefficient of $m^2$ in $x^+(\alpha,m)$
are
\[
\alpha=\alpha^+:=\frac{13284}{16807}\sqrt{2}\simeq1.12,\quad
\alpha=\alpha^-:=-\frac{4428}{2401}\sqrt{2}\simeq-2.61.
\]
Hence, taking for instance $\alpha=\pm {11}/4=\pm 2.75$ we obtain two
values of $\alpha$ for which, for $m$ small enough, the two roots of the
quadratic polynomial in $p_2(x;\alpha,m)$ are positive and so the curve
$C$ is without contact, for $m$ small enough,   in $\{(x,y)\,:\,
x<0\}.$ To know until which values of $m$ this last property holds
we fix $\alpha=\pm 11/4 $ and study the signs of the functions $r_i(\pm
11/4 ,m),$ $i=0,1,2.$

As an example we give some details for
\begin{equation}\label{f1}
r_2(-11/4 ,m)=r_2^0(-11/4 ,m)+r_2^1(-11/4 ,m)\Delta(-11/4 ,m).
\end{equation}
It can be seen that the degrees of $r_2^0(- 11/4 ,m)$ and $r_2^1(-
11/4 ,m)$ are 180 and 184, respectively. Notice that the above expression does
not change sign when
\begin{equation}\label{f2}
\left(r_2^0(- 11/4 ,m)\right)^2-\left(r_2^1(- 11/4
,m)\right)^2(\Delta(- 11/4 ,m))^2
\end{equation}
does not change sign. Hence we have to study the zeroes of this
function which is a polynomial in $m$. This   can be done
analytically by using its Sturm sequence. The smallest positive root
of  the polynomial \eqref{f2} which also vanishes~\eqref{f1} is
$m\simeq 0.40289$. Hence we have proved that for $m\in(0, 4/10),$ it
holds that $r_2(- 11/4 ,m)<0.$ Similarly we can prove that, on the
interval $(0,4/10)$,
\begin{align*}
&r_0( 11/4 ,m)<0,& &r_1( 11/4 ,m)>0,&  &r_2( 11/4 ,m)<0,\\
& r_0( -11/4 ,m)<0,&&r_1(- 11/4 ,m)>0,& &r_2(- 11/4 ,m)<0.
\end{align*}
As a consequence,  the curve $C$  with $b$ given in
\eqref{bfix}, $\alpha=\pm11/4,$ and $m\in(0,4/10)$ is without contact
in $\{(x,y)\,:\, x<0\}.$ It is easy to prove that for $m$
sufficiently small the curve $C$ has a loop in the half-plane
$\{(x,y)\,:\, x<0\}$. Moreover with the same type of algebraic
methods that we have used above it can be seen that the loop also
exists in the interval $m\in (0,4/10)$. Furthermore,  in this
half-plane the loop is crossed transversally  by the flow of
\eqref{btm}, inwards when $\alpha=-11/4$ and outwards when $\alpha=11/4.$
Hence we have proved that when
\begin{equation*}
b=b^{\pm}(m):=\frac5 7 m^2+\frac{72}{2401}m^4\pm\frac{11}4\,m^5,
\end{equation*}
and $m\in(0,4/10)$ it holds that $b^-(m)<b^*(m)<b^+(m).$ From these
inequalities the statement (i) of the  theorem follows by noticing
that $(4/10)^4= 16/625>1/40.$

\subsection{Proof  of Theorem~\ref{mainbt}.(ii)}\label{s3.2} Clearly the ideas of the
method used to prove Theorem~\ref{mainbt}.(i) can be applied when
the curve $C$ given in \eqref{cubica} is taken of higher degree.
Consider now the algebraic curve
\begin{equation}\label{c4}
C_4(x,y):=C_2 (x,y)+\sum_{3\le i+j\le 4} c_{i,j}x^iy^j=0.
\end{equation}
By employing the same procedure applied to the cubic algebraic curve we obtain the following relation between $b$
and $n=m^4$:
\begin{equation}\label{exp}
b=\frac5 7  n^{1/2}+\frac{72}{2401}n- \frac{30024}{45294865}n^{3/2}-
\frac{2352961656}{11108339166925} n^2+O(n^{5/2}).
\end{equation}
 By using the variables $n$ and $b$, the necessary calculations to obtain
 analogous results, for a given
interval $[0,n_0],$ to those obtained in Theorem~\ref{mainbt}.(i)
for the cubic curve are beyond our computational capacity. In
particular we can not obtain a compact expression of the resultant
that would give a proof that $C_4$ is without contact in the
half-plane $\{(x,y)\,:\, x<0\}.$ Nevertheless, once $n$ is fixed we
can perform all  the computations and prove for instance that for
$n=1/20,$
\[\left|\frac5 7  n^{1/2}+\frac{72}{2401}n-
\frac{30024}{45294865}n^{3/2}-b^*(n)\right|<10^{-4}n^{7/4},\] and
for $n=1/100,$
\[
\left|\frac5 7  n^{1/2}+\frac{72}{2401}n-
\frac{30024}{45294865}n^{3/2}- \frac{2352961656}{11108339166925}
n^2-b^*(n)\right|<10^{-5}n^{9/4}.
\]

In any case, we will prove that the expression~\eqref{exp} gives us,
at least locally, the function $b^*(n)$. This can be done by
introducing  the new variables $B>0$ and $M>0$ as
\begin{equation}\label{BM}
B=\frac{b+m^2}2,\quad M^2=\frac{(b+m^2)^2+8m^2}4.
\end{equation}
This change of variables is motivated by equation \eqref{apm}.
Notice that by using them, $a_1^{\pm}=B\pm M.$ Now we repeat all the
procedure developed in the previous subsection but with the curve
$C_4(x,y)=0$ given in \eqref{c4} and these new variables. By forcing
the curve $C_4(x,y)=0$ to coincide as much as possible with both
separatrices of the saddle point (four conditions of contact with
one of them and five with the other one) we determine all the parameters of the quartic.
After that we obtain the relation
\begin{align}\label{bb}
B=&\frac 3 7 M^2-\frac{180} {2401} M^4+ \frac{2366307}{90589730}M^6
-\frac{505643614857}{44433356667700}M^8 +O(M^9),
\end{align}
that ensures one more level of closeness between one of the
separatrices and the algebraic curve, for $M$ sufficiently small.

Fortunately, with these new variables it is now possible to
explicitly get the resultant $P(x;M,B),$ between $C_4(x,y)$ and
\[
D_4(x,y):=\frac{\partial C_4(x,y)}{\partial x}y+\frac{\partial
C_4(x,y)}{\partial y}\left( 2m^2x+(b+m^2)y+x^2+xy\right),
\]
with respect to $y$. We obtain that
\[
P(x;M,B):= \textrm{Res}(C_4,D_4,y)=x^{14}\left(S_0(M,B)\sum_{i=0}^5
T_i(M,B)x^i\right),
\]
where $S_0$ and $T_i,i=0,\ldots,5$ are rational functions in $B$ and
$M.$ To have an idea of their complexity we introduce the following
notation: we will say that a rational function $S(M,B)$ is of type
$\{i,j\}/\{k,l\}$ if after simplifying it, its numerator has
monomials of degree  between $i$ and $j$ and its denominator between
$k$ and $l$. Then
\begin{align*}
&S_0 \quad\mbox{is}\quad \{0,1\}/\{114,132\}, &\\
&T_0 \quad\mbox{is}\quad \{91,120\}/\{0,0\}, \quad &T_1
\quad\mbox{is}\quad
\{109,143\}/\{20,24\},\\
&T_2 \quad\mbox{is}\quad \{124,162\}/\{37,44\}, \quad &T_3
\quad\mbox{is}\quad \{137,179\}/\{52,62\},\\ &T_4
\quad\mbox{is}\quad \{137,179\}/\{53,63\}, \quad &T_5
\quad\mbox{is}\quad \{137,180\}/\{55,65\}.
\end{align*}
Note that  six of the above functions are singular at $B=M=0.$
Inspired by \eqref{bb}  we  express $B$ in terms of $M$ and $\alpha$
as
\begin{align}\label{beta}
B=\beta(M,\alpha):=&\frac 3 7 M^2-\frac{180} {2401} M^4+
\frac{2366307}{90589730}M^6 -\frac{505643614857}{44433356667700}M^8
+\alpha M^9,
\end{align}
where $\alpha$ is an arbitrary  parameter to be determined later. We obtain that
\[
\left.\frac{P(x;\beta(M,\alpha),M)}{x^{14}}\right|_{M=0}=\sum_{i=0}^5
P_i(\alpha)\,x^i,
\]
where now all the functions $P_i,$ $i=0,1,\ldots,5$ are  polynomials
in $\alpha$ with degrees $\{2,5\},\{ 0,6\}, \{0 ,11\},\{0 ,15\},\{0
,15\},$  and we use a  notation similar  to above. We remark that to
have the cancelations between the respective numerators and
denominators in the expression of $P,$ that allow us to evaluate it
at $M=0,$ it is necessary to take $B$ as in the expression
\eqref{beta}.

The polynomial $P_0$ is
\begin{align*}
P_0(\alpha)=k\alpha^2
\big(&296751659628833594552482011388242366232495503625\alpha^3\\
&+719549554938315584569470362390245200816649200\alpha^2\\
&+532161050006783873283311272385459961077760\alpha\\
&+112647235678813465306115059636253208576\big),
\end{align*}
for some positive integer  $k.$ It has $\alpha=0$ as a double root
and three  negative simple roots, which are approximately
$-0.001166$, $-0.000895$ and  $-0.000364$. In particular
\[
P_0(-2500^{-1})<0\quad\mbox{and}\quad P_0({2500}^{-1})>0.
\]
If we fix any of the two values $\alpha^\pm=\pm{2500}^{-1}$ it holds
that $P(x;\beta(M,\alpha^{\pm}),M)$ does not vanish if $M$ and $|x|\ne0$
are small enough. This implies that shrinking, if necessary, the
value of $M$ we obtain that the  oval of the algebraic curve
$C_4(x,y)=0,$ which lies in the negative half plane $\{x<0\}$ and
starts and ends at the origin,  is without contact with the flow of
the system (except at $(0,0)$). Notice that this oval is born from
the cusp point existing when $M=0$ and $B=\beta(\alpha^\pm,0)=0$.

The fact that the resultant between $C_4$ and $D_4$ has different
sign when $\alpha=\alpha^+$ and when $\alpha=\alpha^-$ induces
to think that the flow crosses the oval outwards in one case of plus
sign and inwards in the other one. To prove this fact it suffices to
check this property, for each case, on a single point of the oval.
We have chosen the cutting point of the oval of $C_4(x,y)=0$ and the
$x$-axis.

The four  points of the curve $C_4(x,y)=0$  on the $x$-axis are $0,$
which is a double one, and

\begin{align*}
{\textstyle x_1(\alpha,M)=}&{\textstyle-\frac3 2
M^2+\frac{99}{392}M^4-\frac{14661}{168070}M^6+\frac{1097361567}{28988713600}M^8}\\[2mm]
&{\textstyle-\frac{57336960516777}{3110334966739000}M^{10}+\left(\frac{429}{1280}\alpha
-\frac{1360220314860156764457}{758276772825613000960000}\right)M^{11}+O(M^{12}),}\\[2mm]
{\textstyle x_2(\alpha,M)=}&{\textstyle
-\frac{3266440450}{1212150477}}\\[2mm]
&{\textstyle+\left(-\frac{688506579126746016577}{52895116040087791044}
+\frac{868315765444642387622750}{357042033270592589547}\alpha\right)
M+O(M^2)}.
\end{align*}

Notice that the point on the loop is $(x_1(\alpha,M),0)$. After many
computations we obtain that
\begin{align*}
{\textstyle D_4(x_1(\alpha,M),0)=}&{\textstyle \frac{37182801006000}{184877}\alpha M^{11}}\\
&{\textstyle
-\left(\frac{1233580648740454800}{3107227739}\alpha+\frac{644635919327316538269855546}
{575604427043451909103}\right)M^{12}+O(M^{13})}.
\end{align*}

Hence the sign on $C_4(x,y)=0$ of $D_4$ is different for
$\alpha=\alpha^+>0$ and for $\alpha=\alpha^-<0,$ as we wanted to
prove. Therefore, for $M$ small enough\footnote{To determine an
explicit value of $M$ until which the result holds we could try a
similar study to the one done in the previous subsection, but it
would be extremely long and tedious.}, the value of $B$
corresponding to the existence of the homoclinic loop, say $B^*(M)$,
satisfies
\[
\beta(M,-2500^{-1})<B^*(M)<\beta(M,2500^{-1}).
\]

Collecting all the above results and writing expression~\eqref{bb}
in the old variables $b$ and $n$ we get the proof  of
Theorem~\ref{mainbt}.(ii).

To get an idea of how far the approximation given in
Theorem~\ref{mainbt}.(ii) works  we compare  it  with a numerical
approximation of the bifurcation curve.  Concretely we have obtained
a numerical approximation, working with precision $10^{-16},$ of
$b=b^*(n)$ for $n\in(0,10^{-2})$ by using a Taylor's method for
solving the differential equation. If we denote it by
$b=b^*_{\mbox{\tiny num}}(n)$ it holds that
\begin{align*}
&\left|\frac5 7  n^{1/2}+\frac{72}{2401}n-
\frac{30024}{45294865}n^{3/2}- \frac{2352961656}{11108339166925}
n^2-b^*_{\mbox{\tiny num}}(n)\right|< 6\times 10^{-10},
\end{align*} for
$n\in(0,10^{-2}).$

\subsection{More terms in $b^*(n)$}

By applying again the procedure described in the previous
subsections, but starting with a sextic algebraic curve
\begin{equation*}\label{c5}
C_6(x,y):=C_2 (x,y)+\sum_{3\le i+j\le 6} c_{i,j}x^iy^j=0,
\end{equation*}
we arrive to a relation similar to \eqref{bb} that writes as:
\begin{align*}
B=&{\textstyle \frac 3 7 }M^2-{\textstyle\frac{180} {2401}} M^4+ {\textstyle\frac{2366307}{90589730}}M^6
-{\textstyle\frac{505643614857}{44433356667700}}M^8+{\textstyle\frac{121044460222851597}{21794117111940173000}}M^{10}\\[2mm]
&-{\textstyle\frac{75409774306549331412249}{25960934362572014675870000}}M^{12}+{\textstyle\frac{344552497352535858777709804917}{216470837823465107132247097100000}}M^{14}\\[2mm]
&-{\textstyle\frac{43306460616773431694096161799928995367}{48013140478550999259657196328505023800000}}M^{16}\\[2mm]
&+{\textstyle\frac{91720311301427439675156623493153846098504753619}{174952692015527559148011719369634259959684998000000}}M^{18}\\[2mm]
&-{\textstyle\frac{990748106733217809261982123885784358373281388289053276929}{3187507017704227098311349623537514328886463857480747900000000}}M^{20}\\[2mm]
&+{\textstyle\frac{2174773094339151212022525670857933567647566708078598123403162358137}{11614798115838021630069501971110592530190950393238633306372059000000000}}M^{22}\\[2mm]
&+O(M^{24}).
\end{align*}
From it and the change of variables \eqref{BM} we obtain the
following relation between $b$ and $n$,
\begin{align}
b^*=&{\textstyle\frac5 7 }n^{1/2} +{\textstyle\frac{72}{2401}}n -{\textstyle\frac{30024}{45294865}}n^{3/2}
-{\textstyle\frac{2352961656}{11108339166925}}n^2+{\textstyle\frac{161066618396136}{2724264638992521625}}n^{5/2}\notag\\[2mm]
&-{\textstyle\frac{28575844096870898712}{1622558397660750917241875}}n^3+{\textstyle\frac{37409973403083644863711656}{6764713681983284597882721784375}}n^{7/2}\notag\\[2mm]
&-{\textstyle\frac{1301593321483486009213262204378664}{750205319977359363432143692632890996875}}n^{4}\notag\\[2mm]
&+{\textstyle\frac{750633455019308628819042126726886218707352}{1366817906371309055843841557575267655935039046875}}n^{9/2}\notag\\[2mm]
&-{\textstyle\frac{2188961083333347178341822657596953981848275462851032}{12451199287907137102778709466943415347212749443284171484375}}n^{5}\notag\\[2mm]
&+{\textstyle\frac{1289326941251660725073133052778275691207040442626311930438856}{22685152569996135996229496037325376035529199986794205676507927734375}}n^{11/2}\notag\\[2mm]
&+O(n^{6})\label{conj}.
\end{align}
We believe  that it is an improvement of the expression of $b^*(n)$
given in Theorem~\ref{mainbt}.(ii), with  seven new terms of the
expansion of $b^*(n)$, but we have not been able to prove this fact
due to the complexity of the necessary algebraic calculations.

Notice that the prime decompositions of the denominators  of
\eqref{conj} have a nice and regular  structure:
\begin{align*}
7&=7,\\
 2401&=7^4,\\
45294865&=5\cdot7^7\cdot11,\\
 11108339166925&=5^2\cdot 7^{10}\cdot 11^2\cdot 13,\\
 2724264638992521625&=5^3\cdot 7^{13}\cdot 11^3\cdot 13^2,\\
1622558397660750917241875&=5^4\cdot 7^{15}\cdot 11^4\cdot 13^3\cdot17,\\
67647136819832845978827217843756&=5^5\cdot 7^{18}\cdot 11^5\cdot
13^4\cdot17^2,\\
750205319977359363\cdots92632890996875&=5^5\cdot7^{22}\cdot11^6\cdot13^5\cdot17^3\cdot19,\\
136681790637090558\cdots55935039046875&=5^6\cdot7^{25}\cdot11^7\cdot13^6\cdot17^4\cdot19^2\cdot23,\\
124511992879137102\cdots43284171484375&=5^8\cdot7^{28}\cdot11^8\cdot13^7\cdot17^5\cdot19^3\cdot23^2,\\
226851525699359962\cdots76507927734375&=5^9\cdot7^{31}\cdot11^9\cdot13^8\cdot17^6\cdot19^4\cdot23^3,\\
\end{align*}
where the small irregularity with the number of fives and sevens could be
produced by some cancelations with the respective  numerators.
Unfortunately no regularity appears in the numerators.

This regularity in the denominators of the asymptotic expansion of $b=b^*(n)$ at the origin could give some
clues about a possible closed form expression of this function but, by the moment, we have not been able to obtain it.

\subsection{Relation between our results and Perko's
approach}\label{s3.3}

In \cite{Perko92}, the author  faces the problem of the global
bifurcation diagram, on the Poincar\'{e} sphere, of a different
representation of the Bogdanov-Takens system. Concretely he
considers the system
\begin{equation}\label{sisbt2} \left\{
\begin{array}{ccl}
\dot x&=& y, \\
\dot y&=& x(x-1)+\mu_1y+\mu_2xy
\end{array}
\right.
\end{equation}
and, among other results, he proves the following theorem:

\begin{theorem}[\cite{Perko92}]\label{per} There exists a unique analytic
function $h(\mu_2)$ defined for all $\mu_2\in \mathbb{R}$ that
satisfies $h'(0)=-1/7,$ $h(-\mu_2)=-h(\mu_2),$ and
$\max(-1,-\mu_2^2)<\mu_2h(\mu_2)<0$ for all $\mu_2\ne0$ such that

\begin{enumerate}[(a)]

\item System \eqref{sisbt2} has a unique, hyperbolic limit cycle if and only if $\mu_1
\mu_2<0$ and $0< |\mu_1|<|h (\mu_2)|$; the limit cycle is stable if
$\mu_1> 0,$ and unstable if $\mu_1 < 0$;

\item For $\mu_2\ne 0,$ system \eqref{sisbt2} has a fine focus of multiplicity
one at the origin if and only if $\mu_1 = 0$; for $\mu_2 < 0,$ a
unique, stable, limit cycle is generated in a supercritical Hopf
bifurcation at the origin of \eqref{sisbt2} at the bifurcation value
$\mu_1 =0,$ and it expands monotonically with increasing $\mu_1$
until it intersects the saddle at $(1,0)$ and forms a homoclinic loop
at the bifurcation value $\mu_1 = h(\mu_2)$;

\item System \eqref{sisbt2} has a homoclinic loop at the saddle $(1,0)$ if and
only if $\mu_1 =h(\mu_2)$; the separatrix cycle is hyperbolic if and
only if $\mu_2\ne0$; it is stable (unstable) on its interior for
$\mu_2<0$ ($\mu_2>0$).

\end{enumerate}
\end{theorem}

Let us relate the results of the above theorem with our results. We
need the following simple lemma.

\begin{lemma}\label{lemper} When $\mu_2\ne0$ the change of variables and time
\[
u=\frac{\mu_2^2}{2}\left(2x-1\right),\quad v=\mu_2^3 y, \quad
\tau=\frac{t}{\mu_2}
\]
transforms system \eqref{sisbt2} into
\begin{equation*}\left\{
\begin{array}{ccl}
u'=\dfrac{du}{d\tau}&=& v, \\
v'=\dfrac{dv}{d\tau}&=&-\dfrac{\mu_2^4} 4+
\dfrac{\mu_2(2\mu_1+\mu_2)}{2}v+u^2+uv.
\end{array}
\right.
\end{equation*}
\end{lemma}

As a consequence of the above lemma and Theorems~\ref{teobt},
\ref{mainbt}.(ii) and \ref{per} we obtain:

\begin{corollary}\label{ccc} (i) The following relation holds
\begin{align}
h(\mu_2)&=-\frac{\mu_2}2+\frac{b^*\left(\frac{\mu_2^4}4\right)}{\mu_2}\label{rela}\\&=-\frac1
7\mu_2+\frac{18}{2401}\mu_2^3- \frac{3753}{45294865}\mu_2^5- \frac
{294120207}{22216678333850}\mu_2^7+O(\mu_2^{9}),\label{rela2}
\end{align}
where the functions $b^*$ and $h$ are the ones defined in
Theorems~\ref{teobt} and \ref{per}, respectively.

(ii) The function $b^*(n)$ is given by  $b^*(n)=H(\sqrt{n}),$ for
some analytic function $H$.
\end{corollary}

\begin{proof}
By using Lemma~\ref{lemper} we get that the relations between the
variables in system~\eqref{sisbt2} and the ones of
system~\eqref{sisbt} are
\[
n=\dfrac{\mu_2^4} 4,\qquad b=\dfrac{\mu_2(2\mu_1+\mu_2)}{2}.
\]
By using them we easily obtain the proof, because the relation
$b=b^*(n)$ writes as
\[
\dfrac{\mu_2(2\mu_1+\mu_2)}{2}=b^*\left( \dfrac{\mu_2^4} 4 \right),
\]
which immediately leads to \eqref{rela}. From this relation and the
expression of $b=b^*(n)$ given in Theorem~\ref{mainbt}.(ii), the
expansion~\eqref{rela2} follows. Item~(ii) is a consequence
of~\eqref{rela}, Theorem~\ref{per} and the fact that $h$ is an odd
function.
\end{proof}

We end this section with several remarks.

\begin{remark} (i) Expression \eqref{rela2} improves the local knowledge of the
function~$h$ given in Theorem~\ref{per}. The proof of \cite{Perko92}
that $h'(0)=-1/7$, is equivalent to the classical one which appears
in \cite{GH} and gives the term $\frac5 7 \sqrt{n}$ in
expression~\eqref{exp}. Our approach goes much further and it is
completely different.

(ii) By using Theorem~\ref{mainbt}.(i) we obtain the following
global result: for $|\mu_2|\le 10^{-1/4}\simeq 0.562,$
\[
\left|-\frac 1 7
\mu_2+\frac{18}{2401}\mu_2^3-h(\mu_2)\right|\le\frac{11}{2^{9/2}}\left|\mu_2\right|^4.
\]
Notice that this inequality proves that for
$\mu_2\in(-10^{-1/4},0),$
\begin{equation*}
h(\mu_2)<-\frac17 \mu_2.
 \end{equation*}
This result is coherent with the conjecture made in \cite{Perko92}
that affirms that the above inequality holds for all $\mu_2<0.$

(iii) The inequalities
of Lemma~\ref{btprel} translated to system~\eqref{sisbt2} read as
\[
\max(-1,-\mu_2^2)<\mu_2h(\mu_2)<0.
\]
This information is already contained in Theorem~\ref{per}, but our
proof is different to the one given in \cite{Perko92}.

(iv) In \cite{Perko92} there is another conjecture that says that
$\mu_2h(\mu_2)+1=O(1/\mu_2)$ as $\mu_2\rightarrow-\infty$. In the parameters
of the Bogdanov-Takens system~\eqref{sisbt} it reads as
\[
b^*(n)=\sqrt{n}-1 +O\left({n^{-\frac 1 4}}\right)
\quad\mbox{as}\quad n\rightarrow\infty.
\]

(v) The case $\mu_2=0$ includes new systems which are not contained
in the expression of the Bogdanov-Takens system written as
in~\eqref{sisbt}. For instance the case $\mu_1=\mu_2=0$ corresponds
to a Hamiltonian system with a center at the origin.

\end{remark}

\subsection{A final application}\label{seckorea}

In \cite{korea} it is proved the following result:

\begin{theorem}[\cite{korea}]\label{teok} The system
\begin{equation}\label{sisk} \left\{
\begin{array}{ccl}
\dot x&=& y, \\
\dot y&=& \beta y -\alpha x^2+\alpha^2 x-xy,\quad
\mbox{with}\quad\alpha<0,
\end{array}
\right.
\end{equation}
has a limit cycle if and only if $\gamma<\beta/\alpha<1,$ where
$\gamma$ is a positive constant.
\end{theorem}
This constant $\gamma$ is computed numerically in that paper as
$\gamma\simeq 0.864546.$ We remark that for the values $\alpha$ and
$\beta$ satisfying $\beta/\alpha=\gamma,$ the system~\eqref{sisk}
has a homoclinic loop through the origin. We will improve the
results of that paper. We start by proving the following lemma.

\begin{lemma}\label{lemak}
The constant $\gamma$ defined above is
\[
\gamma=b^*(1/4)+1/2,
\]
where $b^*(n)$ is the function introduced in Theorem~\ref{teobt}.
\end{lemma}

\begin{proof} By applying the change of variables
\[
u=\frac12-\frac x\alpha,\quad v=-\frac y{\alpha^2}, \quad \tau=\alpha
t,
\]
to the system \eqref{sisk} we obtain
\begin{equation*}\left\{
\begin{array}{ccl}
u'=\dfrac{du}{d\tau}&=& v, \\
v'=\dfrac{dv}{d\tau}&=&-\dfrac 14+
\left(\dfrac\beta\alpha-\dfrac12\right)v+u^2+uv.
\end{array}
\right.
\end{equation*}
Hence, the correspondence between the parameters of
system~\eqref{sisk} and the ones of the Bogdanov-Takens
system~\eqref{sisbt} is
\[
n=1/4\quad\mbox{and}\quad b=\beta/\alpha-1/2.
\]
From these relations the lemma follows.
\end{proof}

As we have already explained at the end of Subsection~\ref{s3.2} we
have computed $b^*_{\mbox{\tiny num}}(n)$ at several values of $n$,
obtaining in particular that
\begin{equation}\label{bnum}
\gamma\simeq
    b^*_{\mbox{\tiny num}}(1/4)+1/2=0.864545247421507
=:\gamma_{\mbox{\tiny num}}. \end{equation}
 Notice that only the
first five significative digits agree with that computed
in~\cite{korea}, $0.864546$. As we will prove in
Theorem~\ref{mainkorea}, at least eight of the significative digits
of $\gamma_{\mbox{\tiny num}}$ are correct and $
|\gamma-\gamma_{\mbox{\tiny num}}|<1.2\times 10^{-9}$.

Let us call $\gamma_k:=b^*_k(1/4)+1/2,$ where $b^*_k(n)$ is the
function obtained by adding the first $k$ terms of the asymptotic
expansion of $b^*(n)$ at the origin, see Theorem~\ref{mainbt}.(ii)
and expression~\eqref{conj}. For instance
\[\
b^*_1(n)=\frac57 n^{1/2}, \qquad b^*_2(n)=\frac57
n^{1/2}+\frac{72}{2401}n
\]
and so on. Although we have not proved neither the validity of these
expressions until $n=1/4$ nor the validity of $b^*_k(1/4), k=5,\ldots,11$
we can compute the values $b^*_k(1/4)+1/2$ and compare with the
numerical approximation of $\gamma.$ We obtain
\begin{align*}
{\textstyle\gamma_1= \frac 6 7}&\simeq 0.857142857143,&|\gamma_1-\gamma_{\mbox{\tiny num}}|&\leq 7.5\times 10^{-3}, \\
{\textstyle \gamma_2= \frac{2076}{2401}}&\simeq
0.864639733444,&|\gamma_2-\gamma_{\mbox{\tiny num}}|&\leq 9.5\times 10^{-5},\\
{\textstyle\gamma_3=\frac{39159987}{45294865}} &\simeq 0.864556876370,&|\gamma_3-\gamma_{\mbox{\tiny num}}|&\leq 1.2\times 10^{-5},\\
{\textstyle\gamma_4=\frac{19207287903423}{22216678333850}}&\simeq0.864543637658,&|\gamma_4-\gamma_{\mbox{\tiny num}}|&\leq 1.7\times 10^{-6},\\
{\textstyle\gamma_5=\frac{9421002777077246787}{10897058555970086500}}
&\simeq 0.864545485251,&|\gamma_5-\gamma_{\mbox{\tiny num}}|&\leq
2.4\times
10^{-7},\\
{\textstyle\gamma_6=\frac{11222200726046133491344191}{12980467181286007337935000}}
&\simeq 0.864545210070,&|\gamma_6-\gamma_{\mbox{\tiny num}}|&\leq
3.8\times
10^{-8},\\
{\textstyle\gamma_7=\frac{935744176562\cdots60335616104987}{108235418911\cdots66123548550000}}
&\simeq 0.864545253274,&|\gamma_7-\gamma_{\mbox{\tiny num}}|&\leq
5.9\times
10^{-9},\\
{\textstyle\gamma_8=\frac{2075476618505\cdots38798344748153}{2400657023927\cdots64252511900000}
}&\simeq 0.864545246497,&|\gamma_8-\gamma_{\mbox{\tiny num}}|&\leq 9.2\times 10^{-10},\\
{\textstyle\gamma_9=\frac{7562725921574\cdots99070483500549}{8747634600776\cdots79842499000000}
}&\simeq 0.864545247569,&|\gamma_{9}-\gamma_{\mbox{\tiny num}}|&\leq
1.5\times
10^{-10},\\
{\textstyle\gamma_{10}=\frac{1377872021601\cdots10350610970071}{1593753508852\cdots40373950000000}}&\simeq
0.864545247398,&|\gamma_{10}-\gamma_{\mbox{\tiny num}}|&\leq
2.4\times
10^{-11},\\
{\textstyle\gamma_{11}=\frac{5020759255426\cdots85210508369767}{5807399057919\cdots86029500000000}
}&\simeq 0.864545247425,&|\gamma_{11}-\gamma_{\mbox{\tiny
num}}|&\leq 3.9\times 10^{-12}.
\end{align*}

Notice that the values $\gamma_k$ given above approach well to
$\gamma_{\mbox{\tiny num}},$ and by using again the results of
Theorem~\ref{mainkorea}, also approach to $\gamma.$

In fact in the sequel we will see how the method introduced in this
paper allows to give a concrete interval where the actual value of
$\gamma$ lies.

Recall that our approach for obtaining information of whether a homoclinic connection appears passes trough the construction of the
two polynomials
\begin{align}
C_k(x,y)&=C_2(x,y)+\sum_{3\le i+j\le k} c_{i,j}^kx^iy^j,\label{cubi}\\
D_k(x,y)&=\frac{\partial C_k(x,y)}{\partial x}y+\frac{\partial
C_k(x,y)}{\partial y} (2m^2x+(b+m^2)y+x^2+xy),\notag
\end{align}
with $k\ge3,$ and their resultant. From a computational point of
view it is much simpler the case where all the numbers implied in
their obtention are rational and then $C_k(x,y)$ and $D_k(x,y)$ are
in $\mathbb{Q}[x,y]$.
 It is easy to see that this happens when $b$ and the eigenvalues of the saddle point
 given in \eqref{apm} are rational numbers.
From now one we will particularize our study to
 case $n=1/4,$ that is $m^2=1/2$, although clearly our approach can be adapted to all
 values of $n$ such that $\sqrt{n}$ is a rational number.
In our case, the rationality conditions are reduced to
 \[
b=\frac pq\quad\mbox{and}\quad \left(\frac pq +\frac
12\right)^2+4=\left(\frac r s\right)^2
 \]
for some integer numbers $p,q,r$ and $s.$ This is equivalent to
find integer solutions $p,q$ and $t$ of the quadratic diophantine
equation
\begin{equation}\label{dio}
4p^2+4pq+17q^2=t^2.
\end{equation}
It is well-known how to study this type of equations. More
concretely, if $(p_0,q_0,t_0)\in\mathbb{Z}^3$ is a particular
solution of the diophantine equation
\[
Ap^2+Bpq+Cq^2=Dt^2,
\]
where $A$, $B$, $C$ and $D$ are integers numbers, then
\begin{align*}
p=&p(u,v)=(Ap_0+Bq_0)u^2+2Cq_0uv-Cp_0v^2,\\
q=&q(u,v)=-Aq_0u^2+2Ap_0uv+(Bp_0+Cq_0)v^2,\\
t=&t(u,v)=t_0(Au^2+Buv+Cv^2),
\end{align*}
for any $(u,v)\in\mathbb{Z}^2$, is also an integer solution of the
diophantine equation, because it holds that
\begin{align*}
A\,p(u,v)^2+&B\,p(u,v)\,q(u,v)+C\,q(u,v)^2-D\,t(u,v)^2=\\&(Au^2+Buv+Cv^2)^2
(Ap_0^2+Bp_0q_0+Cq_0^2-Dt_0^2).
\end{align*}

For equation \eqref{dio} it is clear that $(p_0,q_0,t_0)=(1,0,2)$ is
a particular solution. By applying the above procedure we obtain the
new solutions
\begin{align*}
p=&p(u,v)=4u^2-17v^2,\\
q=&q(u,v)=4v(2u+v),\\
t=&t(u,v)=2(4u^2+4uv+17v^2).
\end{align*}
From them we obtain the suitable candidates to perform our study.
We consider, for $(u,v)\in\mathbb{Z}^2,$ with $v(2u+v)\ne0,$
\begin{equation}\label{bbuena}
b=g(u,v)=\frac{4u^2-17v^2}{4v(2u+v)}\in\mathbb{Q},
\end{equation}
for which the eigenvalues of the saddle point given in
\eqref{apm}, $a_1^{\pm},$ are
\[
a_1^+=\frac{2u+v}{4v}\in\mathbb{Q}\quad\mbox{and}\quad a_1^-=-\frac
{4v}{2u+v}\in\mathbb{Q}.
\]

The above results will allow us to prove the following lemma.

\begin{lemma}
Fix $n=1/4$ and $b\in\mathbb{R}.$ Then there are infinitely many
sequences of rational numbers $\{b_j\}_{j\in\mathbb{N}}$ such that
$\lim_{j\to\infty}b_j=b$ and the eigenvalues $a_1^{\pm}$ of the saddle point
\[
\frac{b_j+1/2\pm\sqrt{(b_j+1/2)^2+4}}{2}
\]
are also rational numbers. Moreover these sequences can be
explicitly obtained.
\end{lemma}

\begin{proof} Notice that, for $v\ne0,$ the map $g(u,v)$ given in
\eqref{bbuena} can be written as
\[
b=g(u,v)=\frac{4\left(\frac u v\right)^2-17}{8\left(\frac u v\right)
+4}.
\]

Then the first part of the proof follows by noticing that the graph
of the map
\begin{equation}\label{GG}
w\longmapsto G(w)=\frac{4w^2-17}{8w+4},\quad w\in\mathbb{R}
\end{equation}
covers all the real line. To prove the second part it suffices to
consider any sequence of rational numbers $\{w_j\}_{j\in\mathbb{N}}$
tending to one of the preimages for $G$ of $b,$ say $G^{-1}(b).$
Then, clearly, it holds that for any $j\in\mathbb{N}$
\[b_j:=G(w_j)=\frac{4w_j^2-17}{8w_j+4}\in\mathbb{Q},\quad
\frac{b_j+1/2\pm\sqrt{(b_j+1/2)^2+4}}{2}\in\mathbb{Q}\] and
$\lim_{j\to\infty} b_j=b$ as we wanted to prove.
\end{proof}

Finally, we prove:

\begin{theorem}\label{mainkorea}
Let $\gamma$ be the real number introduced in Theorem~\ref{teok}.
 Then
 \[\gamma\in\left({\textstyle\frac{1127949288}{1304673517},
 \frac{1223980221}{1415750330}}\right)\simeq
 (0.864545247,0.8645452486).\]\end{theorem}

 Notice that the length of the
 interval given in the above theorem is smaller than $1.61\times 10^{-9}.$

\begin{proof} By using Lemma~\ref{lemak}, $\gamma=b^*(1/4)+1/2.$ So
we will fix $n=1/4$ and we study $b^*:=b^*(1/4),$ for
system~\eqref{sisbt}. For the sake of simplicity we only give the
full details of the range of values of $b^*$ that we obtain by
applying our method with $k=3,$ that is by using a cubic
curve~\eqref{cubi}. The proof of the theorem follows by taking
$k=6.$

The numerical approximation \eqref{bnum} of $b^*$ gives  an
orientation for the actual value of $b^*$. After some trials we
consider the values
\[
b^\ell_3:=G\Big(\frac{19}8\Big)=\frac{89}{368}\quad\mbox{and}\quad
b^u_3:=G\Big(\frac{8}3\Big)=\frac{103}{228},
\]
where $G$ is the function given in \eqref{GG}, and we prove that
$b^*\in(b^\ell_3,b^u_3)$. Notice that
$({89}/{368},{103}/{228})\simeq(0.242,0.452)$.

 We will apply the procedure introduced in
Section~\ref{btseccion} for $n=1/4$ and $b\in(b^\ell_3,b^u_3)$ for
searching suitable cubic curves, both having an oval through the
origin and being without contact with the flow of the system.
Moreover we will prove that for one of the ovals the flow goes
inwards and for the other one the flow goes outwards.

Fix one of the values, say $b={89}/{368}$. Then we apply our
procedure to determine $C_3(x,y)$, by imposing all the coincidence
conditions between the cubic and the separatrix associated to the
eigenvalue $a_1^+=23/16.$ Recall that due to our choice of $b$ we
can ensure that all our computations will be with polynomials in
$\mathbb{Q}[x,y].$ We obtain

\begin{align*}
C_3^\ell(x,y)=&{\textstyle\frac{(16x+23y)(16y-23x)}{368}
-\frac{30470974207443747849}{44286028769220680429}x^3
-\frac{48084904789188461109}{88572057538441360858}x^2y}\\
&{\textstyle+\frac{300486284549883520}{44286028769220680429}xy^2
-\frac{78703862917780480}{132858086307662041287}y^3}.
\end{align*}

We omit the explicit expression of
$D^\ell_3(x,y)\in\mathbb{Q}[x,y],$ which has degree 4. The resultant
of $C^\ell_3(x,y)$ and $D^\ell_3(x,y)$, with respect to $y,$
$\textrm{Res}(C^\ell_3,D^\ell_3,y)$ is a polynomial of the form
$x^{10}(Ax^2+Bx+C)\in \mathbb{Q}[x],$ for some huge rational numbers
(numerators and denominators with more than 60 digits), and
$A>0,B<0$ and $C>0$. Hence the resultant is positive in
$\{(x,y)\,:\, x<0\}$.

If we compute the resultant of $C_3^\ell(x,y)$ and $\partial
C_3^\ell(x,y)/\partial y$ with respect to $y$ we obtain a polynomial
of the form $x^2P_4(x),$ with $P_4(x)\in \mathbb{Q}[x].$ The roots
of $P_4$ can be explicitly obtained (or located by using the Sturm
sequences of $P_4$, in view of the application of our approach for
$k>3$). Their approximate values are $-24.478,\, -9.855,\,
x_0\simeq-1.454$ and $32.737.$

The following facts are not difficult to prove for the cubic
$C_3^\ell(x,y)=0$:

\begin{itemize}
\item Its only multiple point is the origin, which is a double point. On
it,
 the curve has two smooth branches tangent to the lines
$(16x+23y)(16y-23x)=0;$

\item The polynomial $C_3^\ell(0,y)$ has the 0 as a double root and a
simple nonzero root;

\item For each $x\in(x_0,0)$ the polynomial $y\longmapsto C_3^\ell(x,y)$ has
exactly three simple real roots;

\item The polynomial $y\longmapsto C_3^\ell(x_0,y)$ has
a simple root and a root of multiplicity two;

\item For each $x\in[-2,x_0)$ the polynomial $y\longmapsto C_3^\ell(x,y)$ has
exactly one simple real root.

\end{itemize}

The tool that we use to prove the first item is the computation of
the resultants between $C_3^\ell(x,y), \partial
C_3^\ell(x,y)/\partial x$ and $\partial C_3^\ell(x,y)/\partial y.$
To prove the other ones we compute the Sturm sequences of the
polynomials $C_3^\ell(x,y)$ considered as polynomials in $y$, with
coefficients in the rational functions with numerators and
denominators in $\mathbb{Q}[x]$. This can be done, except at
finitely many points, given by the zeros of some polynomial of $x$
that appear during the process. In particular,  we can see that
studying the four polynomials $C_3^\ell(x,y)$ for
$x\in\{-2,x_0,-1,0\}$ we have all the information for $y\longmapsto
C_3^\ell(x,y)$ and $x\in[-2,0].$

The above list of properties prove that the curve $C_3^\ell(x,y)=0$
has a loop which starts and ends at the origin being tangent to the
lines $(16x+23y)(16y-23x)=0$ and is contained in the strip
$\{(x,y)\,:\, x_0\le x\le0\}.$ Since we have proved that the
resultant $\textrm{Res}(C^\ell_3,D^\ell_3,y)$ does not vanish on the
left hand plane $\{(x,y)\,:\, x<0\}$ we know that the oval is
without contact for the flow of the system. By studying the sign of
$D^\ell_3$ on the loop we prove that the flow crosses it inwards.
Hence we have proved that $b^*>b^\ell_3={89}/{368}$.

Working similarly, for $b={103}/{228}$ and $a_1={19}/{12}$ we obtain the cubic
\begin{align*}
C_3^u(x,y)=&{\textstyle\frac{(12x+19y)(12y-19x)}{228}
-\frac{42394240475656582327}{66085291294166321043}x^3
-\frac{24891674002318104595}{44056860862777547362}x^2y}\\
&{\textstyle+\frac{320066250082464600}{22028430431388773681}xy^2
-\frac{37723016690931312}{22028430431388773681}y^3}
\end{align*}
and we prove that $b^*<b^u_3={103}/{228}.$

Before considering greater values of  $k$ we want to comment that,
also keeping $k=3,$ but taking rational numbers with big numerators
and denominators, we could improve a little bit the knowledge of the
interval where $b^*$ lies. For instance we can prove that

\[{\textstyle b^*\in\left(G\left(\frac{951}{398}\right),G\left(\frac{29}{11} \right)\right)=
\left(\frac{28898}{114425},\frac{1307}{3036}\right)\simeq
(0.2525,0.4305 ).}\]

By applying our method for $k>3$ we obtain that

\begin{align*}
&{\textstyle
b^*\in\left(\frac{1300991}{3571092},\frac{411011}{1125740}\right)\simeq
(0.364312 ,0.365103)}&\mbox{for}\quad k=4,\\
&{\textstyle
b^*\in\left(\frac{357550843}{980814604},\frac{67268863}{184528084}\right)\simeq(
0.3645448,0.3645454)}&\mbox{for}\quad k=5,\\
&b^*\in\left({\textstyle\frac{951225059}{2609347034},
 \frac{258052528}{707875165}}\right)\simeq
 (0.364545247,0.3645452486)&\mbox{for}\quad k=6.
\end{align*}

We notice that the theorem follows from the result when $k=6.$

All the computations for the six values of $b$ to be studied follow
a similar procedure to the one described above. The main
computational difficulties appear for proving that the algebraic
curve of degree $k$ has a transversal oval passing through the
origin due to  the high degree and huge size of the coefficients of
the polynomials involved. For all cases we prove for the algebraic
curve $C_k(x,y)=0, C_k(x,y)\in\mathbb{Q}[x,y],$ the following
properties, which are similar to the ones given for $k=3.$ In each
of the six cases there is a different value, say $z_0\in(-2,-1),$
which is given as a zero of a polynomial with rational coefficients
computed through a suitable resultant.

\begin{itemize}
\item The only multiple point is the origin, which is a double point. On
it, the curve has two smooth branches tangent to two given lines
with rational slopes;
\item the polynomial $C_k(0,y)$ has the 0 as a double root and $k-2$
simple nonzero roots;
\item for each $x\in(z_0,0)$ the polynomial $y\longmapsto C_k(x,y)$ has
exactly $k$ simple real roots;
\item the polynomial $y\longmapsto C_k(z_0,y)$ has
$k-2$ simple roots and a root of multiplicity two;
\item for each $x\in[-2,z_0)$ the polynomial $y\longmapsto C_k(x,y)$ has
exactly $k-2$ simple real roots.
\end{itemize}
\end{proof}


\begin{thebibliography}{10}

\bibitem{ALGM} A.\,A. Andronov, E.\,A. Leontovich, I.\,I. Gordon and A.\,G. Ma{\u\i}er,
\emph{Theory of bifurcations of dynamic systems on a plane}.
Translated from the Russian. Halsted Press [A division of John Wiley
\& Sons], New York-Toronto, Ont.; Israel Program for Scientific
Translations, Jerusalem-London, 1973.

\bibitem{BFG} B. Boisseau, P. Forg\'{a}cs, and H. Giacomini, \emph{An analytical
approximation scheme to two-point boundary value problems of
ordinary differential equations}. J. Phys. A \textbf{40}
F215--F221, (2007).

\bibitem{CLW} Shui-Nee Chow, Chengzhi Li and Duo Wang, \emph{Normal forms and
bifurcation of planar vector fields}. Cambridge University Press,
Cambridge, 1994.

\bibitem{CGL} B. Coll, A. Gasull and J. Llibre, \emph{Quadratic systems with
a unique finite rest point}. Publ. Mat. \textbf{32}, 199--259,
(1988).

\bibitem{C} W.\,A. Coppel,
\emph{A survey of quadratic systems}. J. Differential Equations
\textbf{2}, 293-304, (1966).

\bibitem{Duff} G.\,F.\,D. Duff, \emph{Limit-cycles and rotated vector fields}. Ann.
of Math. \textbf{67}, 15--31, (1953).

\bibitem{DF} F. Dumortier and P. Fiddelaers, \emph{Quadratic models for generic local
$3$-parameter bifurcations on the plane}. Trans. Amer. Math. Soc.
\textbf{326},  101--126, (1991).



\bibitem{GH} J. Guckenheimer and P. Holmes, \emph{Nonlinear oscillations, dynamical
systems, and bifurcations of vector fields}. Applied Mathematical
Sciences, \textbf{42}. Springer-Verlag, New York, 1983.

\bibitem{korea} Gil-Jun Han, \emph{Bifurcation analysis on an unfolding of
the Takens-Bogdanov singularity}. J. Korean Math. Soc. \textbf{36},
459--467 (1999).


\bibitem{LRW} Chengzhi Li, C. Rousseau and Xian Wang, \emph{Simple proof for
the unicity of the limit cycle in the Bogdanov-Takens system}.
Canad. Math. Bull. \textbf{33},  84--92, (1990).


\bibitem{DLA} F. Dumortier, J. Llibre and J. C.
Art\'{e}s, {\em Qualitative theory of planar differential systems}.
Universitext. Springer-Verlag, Berlin, 2006.


\bibitem{Perko} L.\,M. Perko, \emph{Rotated vector fields and the global behavior of
limit cycles for a class of quadratic systems in the plane}.  J.
Differential Equations  \textbf{18},   63--86, (1975).

\bibitem{Perko92} L.\,M. Perko, \emph{A global analysis of the Bogdanov-Takens system}.
SIAM J. Appl. Math.  \textbf{52}, 1172--1192, (1992).

\bibitem{Per2001}
L.\,M. Perko, \emph{Differential equations and dynamical systems}.
Third ed., Texts in Applied Mathematics, vol.~7, Springer-Verlag,
New York, 2001.


\bibitem{Ro} I.\,G. Rozet, \emph{The closing of the separatrices
of a certain first order differential equation} (Russian)
Differencial'nye Uravnenija  \textbf{7},   2007--2012, (1971).
Translated to English in Differential Equations \textbf{7},
1517--1521, (1971).

\end{thebibliography}
\end{document}